# Поиск равновесий в многостадийных транспортных моделях

*А.В. Гасников (ПреМоЛаб МФТИ, ИППИ РАН)*

*П.Е. Двуреченский (ИППИ РАН, WIAS Berlin)*

*Д.И. Камзолов (ПреМоЛаб МФТИ, ИППИ)*

*Ю.Е. Нестеров (CORE UCL, Belgium)*

*В.Г. Спокойный (ПреМоЛаб МФТИ, ИППИ РАН, ВШЭ, WIAS Berlin)*

*П.И. Стецюк (Институт кибернетики, НАН Украины)*

*А.Л. Суворикова (ПреМоЛаб МФТИ, ИППИ РАН, IRTG 1792 Berlin)*

*А.В. Чернов (Кафедра МОУ ФУПМ МФТИ)*

В работе предлагается оргинальный способ поиска равновесий в многостадийных моделях транспортных потоков. В основе подхода лежит сочетание метода балансировки и универсального градиентного метода. Метод также нашел приложение к задаче поиска барицентра вероятностных мер согласно энтропийно-сглаженному расстоянию Вассерштейна. В последнее время этот круг задач оказался очень популярным в связи с различными приложениями.

**Ключевые слова:** седловая задача, энтропия, задач Монжа–Канторовича, расстояние Вассерштейна, метод баласировки, метод Синхорна, универсальный метод, неточный оракул.

## 1. Введение

Поиск (стохастических) равновесий в многостадийных моделях транспортных потоков приводит к решению следующей седловой задачи с правильной (выпукло-вогнутой) структурой [1 − 4]:

$$\min_{\substack{\sum_{j=1}^{n} x_{ij}=L_i, \sum_{i=1}^{n} x_{ij}=W_j \\ x_{ij}\geq 0,\, i,j=1,\dots,n}} \max_{y\in Q}\left\{\sum_{i,j=1}^{n} x_{ij}\ln x_{ij} + \sum_{i,j=1}^{n} c_{ij}(y)x_{ij} + g(y)\right\} =$$

$$= \max_{y \in Q} \max_{\lambda, \mu \in \mathbb{R}^n} \left\{ \langle \lambda, L \rangle + \langle \mu, W \rangle - \sum_{i,j=1}^{n} \exp\left(-c_{ij}(y) - 1 + \lambda_i + \mu_j\right) + g(y) \right\}, \quad (1)$$

где $g(y)$ и $c_{ij}(y) \geq 0$ – вогнутые гладкие функции (если ищутся не стохастические равновесия, то $c_{ij}(y)$ могут быть негладкими), $Q$ – множество простой структуры, например,

$$Q = \left\{ y : \; y \geq \overline{y} \right\}.$$

Легко понять, что система балансовых ограничений в (1) либо несовместна $\sum_{i=1}^{n} L_i \neq \sum_{j=1}^{n} W_j$, либо вырождена (имеет не полный ранг). В последнем случае это приводит к тому, что двойственные переменные $(\lambda, \mu)$ определены с точностью до произвольной постоянной $C$:

$$\left( \lambda + Ce, \mu - Ce \right), \; e = \underbrace{(1, \ldots, 1)}_{n}.$$

Задачу (1) также можно переписать следующим образом (не ограничивая общности, считаем $\sum_{i=1}^{n} L_i = \sum_{j=1}^{n} W_j = 1$)

$$\min_{\substack{\sum_{j=1}^{n} x_{ij} = L_i, \sum_{i=1}^{n} x_{ij} = W_j \\ x_{ij} \geq 0, i, j = 1, \ldots, n; \; \sum_{i,j=1}^{n} x_{ij} = 1}} \max_{y \in Q} \left\{ \sum_{i,j=1}^{n} x_{ij} \ln x_{ij} + \sum_{i,j=1}^{n} c_{ij}(y) x_{ij} + g(y) \right\} =$$

$$= \max_{y \in Q} \max_{\lambda, \mu \in \mathbb{R}^n} \left\{ \langle \lambda, L \rangle + \langle \mu, W \rangle - \ln\left( \sum_{i,j=1}^{n} \exp\left(-c_{ij}(y) + \lambda_i + \mu_j\right) \right) + g(y) \right\} = -\min_{y \in Q} f(y), \quad (2)$$

где выпуклая функция $f(y)$ определяется как

$$f(y) = \max_{\substack{\sum_{j=1}^{n} x_{ij} = L_i, \sum_{i=1}^{n} x_{ij} = W_j \\ x_{ij} \geq 0, i, j = 1, \ldots, n; \; \sum_{i,j=1}^{n} x_{ij} = 1}} \left\{ -\sum_{i,j=1}^{n} x_{ij} \ln x_{ij} - \sum_{i,j=1}^{n} c_{ij}(y) x_{ij} - g(y) \right\} =$$

$$= \min_{\lambda, \mu} \left\{ \ln\left( \sum_{i,j=1}^{n} \exp\left(-c_{ij}(y) + \lambda_i + \mu_j\right) \right) - \langle \lambda, L \rangle - \langle \mu, W \rangle - g(y) \right\}. \quad (3)$$

Поскольку мы добавили в ограничения условие $\sum_{i,j=1}^{n,n} x_{ij} = 1$, являющееся следствием балансовых уравнений, то это привело к тому, что двойственные переменные $(\lambda, \mu)$ определены с точностью до двух произвольных постоянных $C_\lambda$, $C_\mu$: $(\lambda + C_\lambda e, \mu + C_\mu e)$.

В данной работе мы покажем, как можно решать задачу (2).

Заметим также, что расчет градиента $\nabla f(y)$ (в ряде транспортных приложений вогнутые функции $c_{ij}(y)$ – негладкие, тогда вместо градиентов стоит понимать суперградиенты $c_{ij}(y)$ и субградиент $f(y)$) осуществляется по следующей формуле (Демьянова–Данскина–Рубинова, см., например, [1, 5])

$$\nabla f(y) = -\frac{\sum_{i,j=1}^{n} \exp\left(-c_{ij}(y) + \lambda_i^* + \mu_j^*\right)\nabla c_{ij}(y)}{\sum_{i,j=1}^{n} \exp\left(-c_{ij}(y) + \lambda_i^* + \mu_j^*\right)} - \nabla g(y) = -\sum_{i,j=1}^{n} x_{ij}\left(\lambda^*, \mu^*\right)\nabla c_{ij}(y) - \nabla g(y), \quad (4)$$

где $\left(\lambda^*, \mu^*\right)$ – решение задачи (3), не важно какое именно, градиент $\nabla f(y)$ от выбора $C_\lambda$, $C_\mu$ (см. выше) не зависит. В данной статье мы ограничимся изучением только полноградиентных методов для задачи (2), т.е. не будем рассматривать, например, рандомизацию при вычислении градиента по формуле (4). Планируется отдельно исследовать вопрос о возможности ускорения вычислений за счет введения рандомизации для внешней задачи. На данный момент нам представляется (см. формулу (11) в п. 3), что это может принести дивиденды только в случае, когда вспомогательная задача расчета $\nabla c_{ij}(y)$ достаточно сложная, в свою очередь. Тут требуется много оговорок, в частности, в большинстве приложений умение рассчитывать $\nabla c_{ij}(y)$ для конкретной пары $(i, j)$ без дополнительных затрат позволяет заодно рассчитать и все $\nabla c_{ij}(y)$, $j = 1, ..., n$. Также отдельно планируется исследовать вопрос о том какие подходы и насколько хорошо допускают распараллеливание. Вопрос о целесообразности рандомизации оказывается завязанным и на вопрос о возможности распараллеливания.

Структура статьи следующая. В п. 2 мы рассматриваем популярную в последнее время (в связи с большим числом приложений) задачу вычисления барицентра Вассерштейна различных вероятностных мер. Эта задача оказывается тесно связанной с задачей (1). Мы разбираем в статье этот пример, потому что он хорошо проясняет возможные альтернативы предлагаемому нами основному подходу решения задач (1), (2), изложенному в п. 3. В основе оригинального подхода п. 3 лежит сочетание метода балансировки для решения внутренней задачи оптимизации по двойственным множителям и универсального метода с неточным оракулом для внешней задачи (2). В заключительном п. 4 делаются замечания относительно возможности ускорения подхода, описанного в п. 3.

## 2. Поиск барицентра Вассерштейна

К похожей на (1) задаче приводит поиск барицентра Монжа–Канторовича (в западной литературе чаще говорят барицентра Вассерштейна [6, 7])[1]. Изложим вкратце постановку задачи [7 – 9]. Вводится энтропийно сглаженное транспортное расстояние (см. рис. 1 в [10]), с матрицей $\left\| c_{ij} \right\|_{i,j=1}^{n,n}$, сформированной из квадратов попарных расстояний $c_{ij} = l_{ij}^2$ от носителя меры $i$ до носителя меры $j$ ($L, W \in S_n(1)$):

$$\Delta(L,W) = \min_{\substack{\sum_{j=1}^n x_{ij}=L_i, \sum_{i=1}^n x_{ij}=W_j \\ x_{ij} \geq 0, i,j=1,...,n}} \left\{ \gamma \sum_{i,j=1}^n x_{ij} \ln x_{ij} + \sum_{i,j=1}^n c_{ij} x_{ij} \right\} =$$

$$= \max_{\lambda, \mu} \left\{ \langle \lambda, L \rangle + \langle \mu, W \rangle - \gamma \sum_{i,j=1}^n \exp\left( \frac{-c_{ij} + \lambda_i + \mu_j}{\gamma} - 1 \right) \right\} =$$

$$= \max_{\lambda} \left\{ \langle \lambda, L \rangle \underbrace{- \gamma \sum_{j=1}^n W_j \ln\left( \frac{1}{W_j} \sum_{i=1}^n \exp\left( \frac{-c_{ij} + \lambda_i}{\gamma} \right) \right)}_{H_W^*(\lambda)} \right\}. \tag{5}$$

Определим при $L \in S_n(1)$ функцию $H_W(L) = \Delta(L,W)$. Эта гладкая на $L \in S_n(1)$ функция с градиентом (см. утверждение 3 [9]):

$$\nabla H_W(L) = \lambda^*,$$

где $\lambda^*$ единственное решение (5), удовлетворяющее условию[2] $\langle \lambda^*, e \rangle = 0$. Отсюда следует, что

---

[1] Строго говоря, мы будем искать барицентр вероятностных мер не согласно настоящему (негладкому) расстоянию Вассерштейна (на наш взгляд исторически более правильно это расстояние называть расстоянием Монжа–Канторовича–Добрушина), как можно было подумать из названия, а согласно энтропийно-сглаженному расстоянию Вассерштейна [7].

[2] Решая задачу (6), каким-нибудь прокс-методом с KL прокс-структурой [11], легко понять, что от того, как именно выбирать $\lambda^*$, задаваемое с точностью до сдвига всех компонент на одно и то же произвольное число, метод зависеть не будет. Единственное для чего имеет смысл стремиться к выполнению этого нормирующего условия, так это для лучшей практической обработки (меньшее накопление ошибок округления из-

$$H_W^*(\lambda) = \max_{L \in S_n(1)} \left\{ \langle \lambda, L \rangle - H_W(L) \right\} = \gamma \sum_{j=1}^{n} W_j \ln \left( \frac{1}{W_j} \sum_{i=1}^{n} \exp \left( \frac{-c_{ij} + \lambda_i}{\gamma} \right) \right).$$

Теперь можно перейти к изложению основной конструкции. Задача поиска барицентра Вассерштейна[3] записывается следующим образом:

$$\sum_{k=1}^{m} H_{W_k}(L) \to \min_{L \in S_n(1)}. \tag{6}$$

К сожалению, в такой формулировке мы не можем оценить константу Липшица градиента функционала (6), явно входящую в большинство современных быстрых (ускоренных) численных методов. Однако оказывается (см. п. 3), что существуют быстрые методы, которым для работы не требуется такая информация (константа Липшица градиента).

Перепишем задачу (6), следуя п. 3 работы [9], следующим образом

$$-\sum_{k=1}^{m} H_{W_k}(L_k) \to \max_{\substack{L_1 = L_m | \lambda^1 \\ \cdots\cdots\cdots \\ L_{m-1} = L_m | \lambda^{m-1} \\ L_1, \ldots, L_m \in S_n(1)}},$$

$$\sum_{k=1}^{m-1} \max_{L_k \in S_n(1)} \left\{ \langle \lambda^k, L_k \rangle - H_{W_k}(L_k) \right\} + \max_{L_m \in S_n(1)} \left\{ \left\langle -\sum_{k=1}^{m-1} \lambda^k, L_m \right\rangle - H_{W_m}(L_m) \right\} \to \min_{\lambda^1, \ldots, \lambda^{m-1} \in \mathbb{R}^n},$$

$$\sum_{k=1}^{m-1} H_{W_k}^*(\lambda^k) + H_{W_m}^* \left( -\sum_{k=1}^{m-1} \lambda^k \right) \to \min_{\lambda^1, \ldots, \lambda^{m-1} \in \mathbb{R}^n}, \tag{7}$$

$$L_* = \nabla H_{W_k}^*(\lambda_*^k) \text{ для любого } k = 1, \ldots, m-1,$$

где $L_*$ – единственное решение задачи (6), $\left\{ \lambda_*^k \right\}_{k=1}^{m-1}$ – единственное решение задачи (7). Важное свойство функционала задачи (7) – равномерная ограниченность константы Липшица градиента (следует из [14]). Задача безусловной минимизации (7) может быть эффективно решена различными способами (в зависимости от того насколько велики $n$ и $m$). В частности, для больших $n$ и $m$ неплохо с задачей справляются различные модифи-

---

за конечной длины мантиссы) экспоненциального взвешивания компонент градиента, возникающего на каждом шаге итерационного процесса при выборе KL прокс-структуры.

[3] К сожалению, пока не так много известно о статистической обоснованности использования расстояния Вассерштейна. Другими словами, хотелось бы иметь связь барицентра Вассерштейна с оценками максимального правдоподобия, ну или хотя бы с состоятельными оценками для соответствующих схем экспериментов. Пока установлена только связь с состоятельными оценками [12, 13].

кации метода сопряженных градиентов и быстрых градиентных методов [9]. Структура задачи (7) позволяет эффективно использовать покомпонентные методы (см., например, [15, 16]), которые к тому же хорошо параллелятся для данной задачи. Задача (7) хорошо также решается с помощью распределенных вычислений [17].

В приложениях к поиску разладки требуется много раз перерешивать задачу (7), которую для симметричности перепишем следующим образом

$$\sum_{k=1}^{m} H_{W_k}^{*}\left(\lambda^{k}\right) \to \min_{\substack{\lambda^{1},...,\lambda^{m}\in\mathbb{R}^{n} \\ \sum_{k=1}^{m}\lambda^{k}=0}},$$

немного смещая окно, т.е. заменяя каждый раз несколько первых слагаемых в сумме

$$H_{W_1}^{*}\left(\lambda^{1}\right),...,H_{W_r}^{*}\left(\lambda^{r}\right)$$

на столько же новых (которые, как ожидается, близки к $H_{W_m}^{*}\left(\lambda^{m}\right)$). В таком случае предлагается в итерационном процессе стартовать при сдвиге окошка с того, на чем остановились на прошлом положении окошка, экстраполируя $\lambda_{*}^{m}$ на вновь пришедшие слагаемые. Ясно, что для новой задачи набор

$$\left(\ \lambda_{*}^{r+1},...,\lambda_{*}^{m-1},\lambda_{*}^{m},\underbrace{\lambda_{*}^{m},...,\lambda_{*}^{m}}_{r}\ \right),$$

с которого стартуем, уже не будет оптимальным, однако, мы вправе надеяться на его близость к оптимальному набору, что существенно сокращает число последующих итераций. Интересной, особенно в данном контексте, представляется возможность использования (и интерпретации) распределенных вычислений [17].

Может показаться, что подход, сводящий поиск решения задачи (6) к задаче (7), не доминируем, поскольку в отличие от задачи (6), в задаче (7) мы можем явно выписать функционал и по простым формула рассчитать градиент, который к тому же имеет равномерно ограниченную константу Липшица. С одной стороны, это, действительно, преимущество, но получено оно дорогой ценной – ценной раздутия пространства, в котором происходит оптимизация почти в $m$ раз. И это раздутие скажется не только на сложности одной итерации. Для задачи (6) осуществление одной итерации будет еще более дорогим в виду необходимости на каждой итерации решать $m$ отдельных подзадач расчета $\nabla H_{W_k}\left(L\right)$. Скажется это, прежде всего, на числе необходимых итераций. В следующем разделе будет отмечено, что расчет $\nabla H_{W_k}\left(L\right)$ с помощью метода балансировки не намного сложнее расчета $\nabla H_{W_k}^{*}\left(\lambda^{k}\right)$. При этом задача (6) решается в пространстве намного меньшей размерности, и мы вправе ожидать, что необходимое число итераций может быть на-много меньше, чем для задачи (7). Кроме того, задача (6) решается на компакте, т.е. размер решения (если быть точным, то расстояние от точки старта до решения), входящий в оценку необходимого числа итераций, заведомо ограничен размером симплекса. Задача (7) – задача безусловной оптимизации, причем без свойства сильной выпуклости функционала. Размер ее решения может быть большим, и входит он в оценки необходимого числа итераций также как и для задачи (6), к сожалению, степенным образом (для быстрых (ускоренных) методов можно ожидать линейной зависимости необходимого числа итераций от этого размера). Наконец, для постановок задач об обнаружении разладки (см. выше) также ожидается, что использовать близость решений прямых задач (6) при смещении окошка удастся намного лучше, чем близость в решении двойственных задач (7). В итоге, выгода от подхода, связанного с переходом к задаче (7), уже не столь очевидна, и требует отдельного и более аккуратного исследования, с решающей ролью численных экспериментов.

В ряде задач требуется искать параметрический барицентр Вассерштейна. В таком случае в одном из вариантов постановки предполагают наличие параметрической зависимости $L(\theta) \in S_n(1)$, $\theta \in \Theta$, где размерность вектора параметров $\dim \theta \ll n$. К сожалению, в этом случае нельзя гарантировать с помощью стандартных приемов (стр. 86 [18]) выпуклости задачи

$$\sum_{k=1}^m H_{W_k}\big(L(\theta)\big) \to \min_{\theta \in \Theta}, \qquad (8)$$

за исключением случая, когда $L(\theta) = A\theta + b$, а $\Theta$ – выпуклое множество. В этом случае конструкция (7) видоизменяется следующим образом

$$-\sum_{k=1}^m H_{W_k}\big(L_k\big) \to \max_{\substack{L_1 = L_m |\lambda^1 \\ \cdots\cdots\cdots \\ L_{m-1} = L_m |\lambda^{m-1} \\ L_m = A\theta + b |\tilde{\lambda} \\ L_1, \ldots, L_m \in S_n(1), \theta \in \Theta}},$$

$$\sum_{k=1}^{m-1} H_{W_k}^*\big(\lambda^k\big) + H_{W_m}^*\left(-\sum_{k=1}^{m-1} \lambda^k - \tilde{\lambda}\right) + \big\langle \tilde{\lambda}, A\theta + b \big\rangle \to \min_{\substack{\lambda^1, \ldots, \lambda^{m-1}, \tilde{\lambda} \in \mathbb{R}^n \\ \theta \in \Theta}}, \qquad (9)$$

$$L_* = \nabla H_{W_k}^*\big(\lambda_*^k\big) \text{ для любого } k = 1, \ldots, m-1.$$

Как следствие, нет никаких гарантий, что изложенная выше конструкция, связанная с переходом к двойственной задаче (7), и восстановлению решения прямой задачи (6) по явным формулам через двойственные множители, в общем случае будет работать хотя бы для поиска локальных решений.[4] Другими словами, необходимо искать глобальный оптимум задачи (8), исходя из работы с прямой задачей (8). Один из вариантов того, как это можно делать, будет описан в следующем пункте.[5]

Однако при другом варианте постановки (более предпочтительном) можно задавать зависимость $L(\theta)$ с помощью аффинных равенств и выпуклых неравенств

$$\sum_{k=1}^{m} H_{W_k}\left(L\right) \to \min_{\substack{A\theta+BL=c \\ g(\theta,L)\leq0 \\ L\in S_n(1); \theta\in\Theta}} .$$

Многие параметрические зависимости можно загнать в такое представление [20]. В частности, отметим полиэдральные представления Фурье–Моцкина [20], возникающие, например, в робастной оптимизации

$$\sum_{k=1}^{m} H_{W_k}\left(L\right) \to \min_{L\in\{L\in S_n(1):\exists\theta: A\theta+BL\leq c\}} .$$

Можно переписать задачи (7), (9) и на эти случаи, причем сделать это корректно в том смысле, что правомочность подхода полностью сохранится. При этом принципиально ничего из сказанного выше не поменяется. Подробнее об этом планируется написать в отдельной работе.

В действительности, в приложениях наиболее интересен случай, когда ищется барицентр именно расстояний Вассерштейна,[6] а не энтропийно-сглаженных расстояний [7 – 10]. Другими словами, интересно изучать предельное поведение $\gamma \to 0+$ (см. п. 3.1 [8], утверждение 1 [9], п. 3 и конец п. 4 [10]). К сожалению, методы из пп. 2, 3 оказываются весьма чувствительными к этому предельному переходу. Для метода из этого раздела кон-станта Липшица градиента в задаче (7) будет расти как $\gamma^{-1}$, соответственно, число итераций будет увеличиваться (при использовании быстрых (ускоренных) методов) как $\gamma^{-1/2}$. Еще более плохое поведение (см. [5]) можно ожидать от метода балансировки, использующегося в подходе из п. 3. Планируется в отдельной публикации исследовать вопрос о том, как следует действовать при малых $\gamma > 0$. В частности, в вырожденном случае $\gamma = 0$. По-видимому, в этом случае поможет философия искусственного сглаживания[7] [14], в которой искусственно введенная энтропийная регуляризация уже задается с четко заданным коэффициентом регуляризации $\gamma > 0$, зависящим от итоговой точности по функции, с которой требуется решить задачу. Другой способ – использовать менее чувствительные (чем метод балансировки) способы решения двойственной задачи, например, при небольших значениях $n$ ожидается, что лучше сработает $r$-метод Н.З. Шора и некоторые его обобщения [23, 24]. В данной работе мы фиксируем $\gamma > 0$ и далее уже не будем возвращаться к подобного рода вопросам.

---

[4] Впрочем, есть результаты (см. формулу (8) п. 3 § 2 главы 8 [19]) о локальной сходимости обычного градиентного спуска для задачи (9) при некоторых дополнительных предположениях.

[5] При этом правомочность подхода п. 3 для постановки задачи (8) имеет место при дополнительном предположении, что метод стартует из выпуклой окрестности точки минимума с небольшим запасом, допускающим возможность по ходу итерационного процесса оказаться дальше от решения, чем в начальный момент.

[6] Численные методы поиска "честного" барицентра Вассерштейна вероятностных мер в основном строятся на том, что когда меры заданы на прямой, существуют эффективные способы решения задачи поиска барицентра [6]. Далее проектируют (считают преобразования Радона) меры на случайные прямые и решают одномерные задачи. По их решениям восстанавливают решение исходной задачи [21, 22]. В отличие от других подходов, здесь существенно используется структура матрицы затрат $c_{ij}=l_{ij}^2$ (в дискретном случае). Интересно было бы исследовать вопрос о применимости этого подхода к постановкам задач о разладках, в которых требуется много раз пересчитывать барицентр (см. выше). Также интересно было бы сравнить описанные подходы с остальными. Этому планируется посвятить отдельную публикацию.

В заключение этого раздела отметим, что поиск барицентра Вассерштейна в случае $m = 1$ может быть осуществлен явно: $L = W$. Обоснование этого частного результата представляется довольно полезным для понимания основной конструкции этого раздела.

### 3. Универсальный метод с неточным оракулом

Из п. 2 следует, что внутренняя задача максимизации по $(\lambda, \mu)$ может быть явно решена по $\mu$ при фиксированном $\lambda$, и наоборот (это верно для задач (1) и (2), и приводит к одним и тем же формулам). Собственно, таким образом, получается метод балансировки расчета матрицы корреспонденций по энтропийной модели, см., например, [5] (тесно связанный с методом Синхорна [8, 10]), как метод простой итерации для явно выписываемых условий экстремума (принципа Ферма): $\lambda = \Lambda(\mu)$, $\mu = \mathrm{M}(\lambda)$.

---

[7] Это сглаживание правильно называть двойственным сглаживанием, поскольку для того чтобы добиться гладкости в прямой негладкой задаче, которая имеет Лежандровское (седловое) представление [14, 20], в это представление, которое также можно понимать как двойственное, вводят аддитивным образом с небольшим коэффициентом сильно выпуклый (вогнутый) функционал. Этот функционал и обеспечивает гладкость (а еще точнее Липшицевость градиента) в прямой задаче. В нашем случае мы исходим из задачи о перемещении масс (Монж–Канторович), являющейся задачей ЛП. Однако для удобства вычисления расстояний Вассерштейна мы перешли к двойственной задаче. Мы хотим сделать гладкой двойственную задачу, потому что именно с ней в дальнейшем и идет работа. Для этого двойственное сглаживание (в нашем случае энтропийное) применяется к двойственной задаче для двойственной задачи к транспортной задаче, т.е. применяется просто к транспортной задаче.

Метод балансировки имеет вид ($\left[\lambda\right]_0 = \left[\mu\right]_0 = 0$):

$$\left[\lambda_i\right]_{k+1} = -\ln\left(\frac{1}{L_i}\sum_{j=1}^{n}\exp\left(-c_{ij}\left(y\right)-1+\left[\mu_j\right]_k\right)\right),$$

$$\left[\mu_j\right]_{k+1} = -\ln\left(\frac{1}{W_j}\sum_{i=1}^{n}\exp\left(-c_{ij}\left(y\right)-1+\left[\lambda_i\right]_k\right)\right)$$

или

$$\left[\mu_j\right]_{k+1} = -\ln\left(\frac{1}{W_j}\sum_{i=1}^{n}\exp\left(-c_{ij}-1+\left[\lambda_i\right]_{k+1}\right)\right).$$

В этих формулах "–1" в экспоненте для метода (2) (в отличие от метода (1)) можно не писать, поскольку двойственные множители определяются неоднозначным образом с бо́льшим произволом для задачи (2) (см. выше), достаточным для справедливости этого замечания.

Оператор $\left(\lambda,\mu\right)\rightarrow\left(\Lambda\left(\mu\right),\mathrm{M}\left(\lambda\right)\right)$ является сжимающим в метрике Биркгофа–Гильберта $\rho$ [25]. Это означает, после $N\sim\ln\left(\sigma^{-1}\right)$ итераций метода балансировки можно получить такие $\left(\lambda_N,\mu_N\right)$, что ($\left\{\left(\lambda_*\left(y\right),\mu_*\left(y\right)\right)\right\}$ – двумерное аффинное множество решений, см. п. 1)

$$\rho\left(\left(\lambda_N,\mu_N\right);\left\{\left(\lambda_*\left(y\right),\mu_*\left(y\right)\right)\right\}\right)\leq\sigma. \tag{10}$$

Причем на практике наблюдается очень быстрая сходимость, т.е. коэффициент пропорциональности не большой [5]. Таким образом, мы можем приближенно решить внутреннюю задачу.[8]

Далее предлагается воспользоваться прямо-двойственным (эта важна, поскольку нужно восстанавливать двойственные переменные) универсальным методом [26] для ре-шения внешней задачи оптимизации по $y$ (имеется видео/презентация с описанием этого метода [27]). К сожалению, в формулировке (1) (в отличие от формулировки (2)) кроме того что внешняя задача гладкая (при условии гладкости $c_{ij}(y)$ [28, 29]), больше ничего о ней сказать нельзя (константа Липшица градиента не ограничена). Также не понятна гладкость задачи (6). Поэтому и по ряду других причин, о которых будет сказано далее, было отдано предпочтение универсальному методу, оптимально адаптивно настраивающемуся на гладкость функционала $f(y)$ на текущем участке пребывания итерационного процесса.[9] Однако нам потребуется использовать этот метод в варианте с неточным оракулом, выдающим градиент [30].

---

[8] Заметим, что в пп. 3.1, 3.2 работы [8] предлагается за счет раздутия прямого пространства с помощью обобщения описанного метода балансировки Брэгмана (метода проекций Брэгмана) решать задачу поиска барицентра напрямую, т.е. отпадает необходимость в решении внешней задачи. Плата за это достаточно большая – увеличение размера прямого пространства в $m$ раз, но метод при этом хорошо параллелится.

Напомним (см. п. 1), что мы решаем задачу (2), представимую в виде (здесь в max представлении $x = x$, $\bar{Q} = \left\{ x_{ij} \geq 0, \ i, j = 1, ..., n : \sum_{j=1}^{n} x_{ij} = L_i, \sum_{i=1}^{n} x_{ij} = W_j \right\}$, а в min представлении $x = (\lambda, \mu)$, $\bar{Q} = \mathbb{R}^{2n}$ — см. формулу (3)):

$$f(y) = \max_{x \in \bar{Q}} \Psi(x, y) = \min_{x \in Q} \Phi(x, y) \to \min_{y \in Q}.$$

Далее везде будем предполагать, что $y \in Q$.

**Определение 1 (см. главу 4 [31]).** $(\delta, L)$-*оракул выдает (на запрос, в котором указывается только одна точка $y$) такие $\left( F(y), G(y) \right)$, что и для любых $y, y' \in Q$*

$$0 \leq f(y') - F(y) - \langle G(y), y' - y \rangle \leq \frac{L}{2} \|y' - y\|^2 + \delta.$$

---

[9] Бытует мнение, с которым столкнулись и авторы данной статьи, что любой универсальный метод должен чем-то платить за свою универсальность, и в этой связи возникает много вопросов, в частности: насколько дорога эта плата? В принципе, в статье [26] довольно подробно проясняется этот момент. Тем не менее, мы повторим здесь соображения из [26]. Действительно, плата за универсальность есть. Универсальный метод из работы [26] может сделать где-то в 4 раза больше итераций для задачи с, более, менее, одинаковой константой Липшица градиента во всей области (где довелось пройти итерационному процессу), по сравнению с обычным быстрым градиентным методом [14]. Тем не менее, замечательная особенность универсального метода не только в том, что он настраивается на гладкость задачи и применимым к любым задачам, но и в том, что (в отличие от подавляющего большинства методов) этот метод локально настраивается на гладкость функционала. И для сильно неоднородных функционалов типично, что универсальный метод делает заметно меньше итераций, чем, скажем, быстрый градиентный метод (плата за это уже учтена в отмеченном выше потенциально возможном увеличении числа итераций в 4 раза в худшем случае). Примеры, поясняющие сказанное, имеются в работе [26].

Из определения 1 разу следует, что для любого $x \in Q$

$$F(x) \le f(x) \le F(x) + \delta$$

и для любых $x, y \in Q$

$$f(y) \ge f(y) - \langle G(x), y - x \rangle - \delta .$$

Из последнего свойства получаем, что определение $(\delta, L)$-оракула можно понимать как обобщение на гладкие задачи классического понятия негладкой выпуклой оптимизации: $\delta$-субградиента (см. п. 5 § 1 главы 5 [19]). В приводимом далее утверждении в первой его части следует сохранить обозначения для задачи (2), (3) и следует обозначить $x = \lambda$, $y = L$ для задачи (5), (6); а во второй части утверждения следует обозначить $x = (\lambda, \mu)$, $y = y$ для задачи (2), (3). Таким образом, на задачу (2), (3) можно посмотреть с двух разных ракурсов, однако второй ракурс менее привлекателен в виду необходимости рассмотрения ограниченых множеств $\bar{Q}$, что в интересующих нас приложениях место не имеет.

**Утверждение 1.** *Если* $\psi(y) = \max\limits_{x \in \bar{Q}} \Psi(x, y)$, *где* $\Psi(x, y)$ – *выпуклая по* $y$ *и вогнутая по* $x$ *функция, и найден такой* $\tilde{x} \in \bar{Q}$, *что*

$$\psi(y) - \Psi(\tilde{x}, y) \le \delta ,$$

*то субградиент* $\partial_y \Psi(\tilde{x}, y)$ – *есть* $\delta$-*субградиент функции* $\psi(y)$ *в точке* $y$.

*Если* $\varphi(y) = \min\limits_{x \in Q} \Phi(x, y)$, *где* $\Phi(x, y)$ – *выпуклая по совокупности переменных функция, и найден такой* $\tilde{x} \in \bar{Q}$, *что*

$$\max\limits_{z \in Q} \langle \Phi_x(\tilde{x}, y), \tilde{x} - z \rangle \le \delta ,$$

*то*

$$\Phi(\tilde{x}, y) - \varphi(y) \le \delta$$

*и субградиент* $\Phi_y(\tilde{x}, y) = \partial_y \Phi(\tilde{x}, y)$ – *есть* $\delta$-*субградиент функции* $\varphi(y)$ *в точке* $y$.

**Доказательство.** Ограничимся доказательством только второй части этого утверждения. Доказательство первой части см. на стр. 124 (лемма 13) книги [19].

Из выпуклости $\Phi(x, y)$ по совокупности переменных имеем

$$\Phi(x', y') \geq \Phi(x, y) + \langle \Phi_x(x, y), x' - x \rangle + \langle \Phi_y(x, y), y' - y \rangle. \qquad (*)$$

Определим зависимость $x(y)$ из соотношения

$$\varphi(y) = \min_{x \in Q} \Phi(x, y) = \Phi(x(y), y).$$

Заметим, что $\Phi(\tilde{x}, y) \geq \varphi(y)$. Положим в (*) $x' = x(y')$, $x = \tilde{x}$. Тогда

$$\varphi(y') = \Phi(x', y') \geq \Phi(\tilde{x}, y) + \langle \Phi_x(\tilde{x}, y), x' - \tilde{x} \rangle + \langle \Phi_y(\tilde{x}, y), y' - y \rangle \geq$$

$$\geq \varphi(y) + \langle \Phi_x(\tilde{x}, y), x(y') - \tilde{x} \rangle + \langle \Phi_y(\tilde{x}, y), y' - y \rangle \geq \varphi(y) + \langle \Phi_y(\tilde{x}, y), y' - y \rangle - \delta.$$

В последней формуле мы использовали, что

$$\langle \Phi_x(\tilde{x}, y), \tilde{x} - x(y') \rangle \leq \delta.$$

В свою очередь, из выпуклости $\Phi(x, y)$ по $x$ (для всех допустимых $y$), имеем

$$\Phi(\tilde{x}, y) - \Phi(x(y'), y) \leq \langle \Phi_x(\tilde{x}, y), \tilde{x} - x(y') \rangle.$$

Беря в этой формуле $y' = y$ и воспользовавшись определением $x(y)$, получаем, что

$$\Phi(\tilde{x}, y) - \varphi(y) \leq \langle \Phi_x(\tilde{x}, y), \tilde{x} - x(y) \rangle. \; \bullet$$

Однако не хочется довольствоваться возможностью находить только $\delta$-субградиент (из утверждения 1 эта возможность очевидна), поскольку в определенных ситуациях явно можно рассчитывать на некоторую гладкость итоговой (внешней) задачи (2). Понятие $(\delta, L)$-оракула в некотором смысле налагает наиболее слабые условия на возможные неточности в вычислении функции и градиента, при которых можно рассчитывать, что скорость сходимости метода, учитывающего гладкость (Липшицевость градиента функционала) задачи, сильно не пострадает (см. теорему 1 ниже).

На первый взгляд может показаться, что применимость описанной концепции $(\delta, L)$-оракула к задаче (1) следует из следующего результата (см. п. 4.2.2 [31]).

**Утверждение 2.** *Пусть подзадача энтропийно-линейного программирования (ЭЛП) в (2) решена (по функции) с точностью $\delta$, т.е. найден такой $\tilde{x}(c)$, удовлетворяющей балансовым ограничениям, что*

$$\sum_{i,j=1}^{n} \tilde{x}_{ij}(c) \ln \tilde{x}_{ij}(c) + \sum_{i,j=1}^{n} c_{ij} \tilde{x}_{ij}(c) - \min_{\substack{\sum_{j=1}^{n} x_{ij}=L_i, \sum_{i=1}^{n} x_{ij}=W_j \\ i,j=1,\ldots,n}} \left\{ \sum_{i,j=1}^{n} x_{ij} \ln x_{ij} + \sum_{i,j=1}^{n} c_{ij} x_{ij} \right\} \leq \delta.$$

*Тогда для функции*

$$\bar{f}(c) = -\min_{\substack{\sum_{j=1}^{n} x_{ij}=L_i, \sum_{i=1}^{n} x_{ij}=W_j \\ i,j=1,\ldots,n}} \left\{ \sum_{i,j=1}^{n} x_{ij} \ln x_{ij} + \sum_{i,j=1}^{n} c_{ij} x_{ij} \right\}$$

*набор*

$$-\left( \sum_{i,j=1}^{n} \tilde{x}_{ij}(c) \ln \tilde{x}_{ij}(c) + \sum_{i,j=1}^{n} c_{ij} \tilde{x}_{ij}(c), \left\{ \tilde{x}_{ij}(c) \right\}_{i,j=1}^{n,n} \right)$$

*является $\left( \delta, 2 \cdot \max\limits_{i,j=1,\ldots,n} c_{ij} \right)$-оракулом.*

К сожалению, большинство методов (в том числе метод балансировки) не удовлетворяют одному пункту утверждения 2, а именно, они выдают вектор $\tilde{x}$, который лишь приближенно удовлетворяет балансовым ограничениям (в утверждении требование точного удовлетворения балансовых ограничений является существенным, и не может быть как-то равнозначно релаксировано). Связанно это с тем, что для задачи ЭЛП, когда ограничений намного меньше числа прямых переменных, обычно решается двойственная задача, по которой восстанавливается решение прямой задачи [5, 32]. Как следствие приобретается невязка и в ограничениях. Собственно, в представлении градиента функционала по формуле (4) имеются два способа. Первый через двойственные множители $(\lambda, \mu)$, второй через решение прямой задачи $x$. Функционал прямой задачи сильно выпуклый по $x$, поскольку энтропия 1-сильно выпуклая функция в 1-норме [14]. Поэтому сходимость в решении прямой задачи по функции обеспечивает сходимость и по аргументу, что и означает возможность определения с хорошей точностью градиента по формуле (4) через $x$. Другая ситуация возникает, если смотреть на двойственную задачу к задаче ЭЛП (в приводимом далее утверждении следует обозначить $x = (\lambda, \mu)$, $y = y$ для задачи (2), (3)).

**Утверждение 3.** *Пусть $\varphi(y) = \min\limits_{x \in Q} \Phi(x, y)$, где $\Phi(x, y)$ – такая достаточно гладкая, выпуклая по совокупности переменных функция, что*[10]

---

[10] Это утверждение имеет достаточно простую геометрическую интерпретацию. Проекция надграфика выпуклой функции будет выпуклым множеством, то есть, в свою очередь, надграфиком некоторой выпуклой функции. Кривизна границы у полученного при проектировании множества будет не больше, чем была у исходного множества. Это следует из того, что проектирование – сжимающий оператор.

$$\left\| \nabla \Phi \left( x', y' \right) - \nabla \Phi \left( x, y \right) \right\|_2 \le L \left\| \left( x', y' \right) - \left( x, y \right) \right\|_2.$$

*Пусть для произвольного $y \in Q$ (считаем, что множество $Q$ содержит внутри себя шар радиуса более $\sqrt{2\delta/L}$) можно найти такой $\tilde{x} = \tilde{x}(y) \in \bar{Q}$, что*

$$\max_{z \in \bar{Q}} \left\langle \Phi_x \left( \tilde{x}, y \right), \tilde{x} - z \right\rangle \le \delta.$$

*Тогда*

$$\Phi \left( \tilde{x}, y \right) - \varphi \left( y \right) \le \delta,$$

$$\left\| \nabla \varphi \left( y' \right) - \nabla \varphi \left( y \right) \right\|_2 \le L \left\| y' - y \right\|_2,$$

*и $\left( \Phi \left( \tilde{x}, y \right) - 2\delta, \Phi_y \left( \tilde{x}, y \right) \right)$ будет $\left( 6\delta, 2L \right)$-оракулом для $\varphi \left( y \right)$ на выпуклом множестве, полученном из множества $Q$ отступанием от границы $\partial Q$ во внутрь $Q$ на расстояние $\sqrt{2\delta/L}$ (по условию это множество не пусто).*

**Доказательство.** По условию задачи имеем при всех допустимых значениях аргументов $\Phi$:

$$\lambda_{\max} \left( \left\| \begin{matrix} \Phi_{xx} & \Phi_{xy} \\ \Phi_{yx} & \Phi_{yy} \end{matrix} \right\| \right) = \sup_{\|h\|_2 \le 1} \left\langle h, \left\| \begin{matrix} \Phi_{xx} & \Phi_{xy} \\ \Phi_{yx} & \Phi_{yy} \end{matrix} \right\| h \right\rangle \le L. \tag{**}$$

Заметим, что также по условию при всех допустимых значениях аргументов $\Phi$:

$$\left\| \begin{matrix} \Phi_{xx} & \Phi_{xy} \\ \Phi_{yx} & \Phi_{yy} \end{matrix} \right\| \succ 0, \ \Phi_{xx} \succ 0, \ \Phi_{yy} \succ 0, \ \Phi_{yx} = \Phi_{xy}^T, \ \Phi_{xx} = \Phi_{xx}^T, \ \Phi_{yy} = \Phi_{yy}^T.$$

Для упрощения последующих рассуждений (в частности, чтобы не работать с псевдообратными матрицами) будем, считать, что матрица $\Phi_{xx} \succ 0$ положительно определена (исходя из условий, гарантировать можно лишь неотрицательную определенность). Также будем считать (в интересующих нас приложениях к задачам (2), (6) это имеет место), что зависимость $x(y)$, определяемая из соотношения

$$\varphi \left( y \right) = \min_{x \in Q} \Phi \left( x, y \right) = \Phi \left( x \left( y \right), y \right)$$

однозначным образом, и удовлетворяет соотношению

$$\Phi_x \left( x \left( y \right), y \right) \underset{y}{\equiv} 0,$$

из которого имеем

$$\Phi_{xx}\big(x(y),y\big)\left\|\frac{\partial x(y)}{\partial y}\right\|+\Phi_{xy}\big(x(y),y\big)\underset{y}{\equiv}0,$$

т.е.

$$\|\partial x/\partial y\|=\|\partial x_i/\partial y_j\|=-\Phi_{xx}^{-1}\Phi_{xy}.$$

Поскольку $\varphi(y)=\Phi\big(x(y),y\big)$, то

$$\varphi_{yy}=\|\partial x/\partial y\|^T\Phi_{xx}\|\partial x/\partial y\|+\|\partial x/\partial y\|^T\Phi_{xy}+\Phi_{yx}\|\partial x/\partial y\|+\Phi_{yy}=\Phi_{yy}-\Phi_{yx}\Phi_{xx}^{-1}\Phi_{xy}.$$

С учетом этой формулы и из формулы дополнения по Шуру [33], получаем

$$\left\|\begin{matrix}\Phi_{xx}&\Phi_{xy}\\\Phi_{yx}&\Phi_{yy}\end{matrix}\right\|=\left\|\begin{matrix}E_x&0\\\Phi_{yx}\Phi_{xx}^{-1}&E_y\end{matrix}\right\|\left\|\begin{matrix}\Phi_{xx}&0\\0&\varphi_{yy}\end{matrix}\right\|\left\|\begin{matrix}E_x&\Phi_{xx}^{-1}\Phi_{xy}\\0&E_y\end{matrix}\right\|,$$

где $E_x$, $E_y$ – единичные матрицы соответствующих размеров. Поскольку

$$\left\|\begin{matrix}E_x&0\\\Phi_{yx}\Phi_{xx}^{-1}&E_y\end{matrix}\right\|=\left\|\begin{matrix}E_x&\Phi_{xx}^{-1}\Phi_{xy}\\0&E_y\end{matrix}\right\|^T,$$

и эти матрицы полного ранга, то из (**) имеем, что

$$\sup_{\|h\|_2\leq1}\big\langle h,\varphi_{yy}\,h\big\rangle=\lambda_{\max}\big(\Phi_{yy}\big)\leq\max\big\{\lambda_{\max}\big(\Phi_{xx}\big),\lambda_{\max}\big(\varphi_{yy}\big)\big\}=\lambda_{\max}\left(\left\|\begin{matrix}\Phi_{xx}&\Phi_{xy}\\\Phi_{yx}&\Phi_{yy}\end{matrix}\right\|\right)\leq L.$$

Таким образом, установлено, что

$$\varphi(y)\leq\varphi(x)+\big\langle\nabla\varphi(x),y-x\big\rangle+\frac{L}{2}\|y-x\|_2^2.$$

Согласно утверждению 1

$$\varphi(y)\geq\varphi(x)+\big\langle\Phi_y(\tilde{x},y),y-x\big\rangle-\delta.$$

Далее проведем рассуждения аналогично стр. 107 (и не много отлично от стр. 115) диссертации [31]. Вычитая из первого неравенства второе, получим

$$\big\langle\Phi_y(\tilde{x},y)-\nabla\varphi(x),y-x\big\rangle\leq\frac{L}{2}\|y-x\|_2^2+\delta.$$

Положим ( $t > 0$ )

$$y - x = t \frac{\Phi_y(\tilde{x}, y) - \nabla \varphi(x)}{\left\| \Phi_y(\tilde{x}, y) - \nabla \varphi(x) \right\|_2},$$

получим

$$\left\| \Phi_y(\tilde{x}, y) - \nabla \varphi(x) \right\|_2 \le \frac{Lt}{2} + \frac{\delta}{t}.$$

Минимизируя правую часть неравенства по $t > 0$, получим (при $t = \sqrt{2\delta/L}$), что

$$\left\| \Phi_y(\tilde{x}, y) - \nabla \varphi(x) \right\|_2 \le \sqrt{2\delta L}.$$

Отсюда и из утверждения 1 имеем, что

$$\varphi(y) \le \varphi(x) + \left\langle \nabla \varphi(x), y - x \right\rangle + \frac{L}{2} \| y - x \|_2^2 \le$$

$$\le \varphi(x) + \left\langle \Phi_y(\tilde{x}, y), y - x \right\rangle + \sqrt{2\delta L} \| y - x \|_2 + \frac{L}{2} \| y - x \|_2^2 \le$$

$$\le \Phi(\tilde{x}, y) - 2\delta + \left\langle \Phi_y(\tilde{x}, y), y - x \right\rangle + \sqrt{2\delta L} \| y - x \|_2 + \frac{L}{2} \| y - x \|_2^2 + 2\delta \le$$

$$\le \Phi(\tilde{x}, y) - 2\delta + \left\langle \Phi_y(\tilde{x}, y), y - x \right\rangle + \frac{2L}{2} \| y - x \|_2^2 + 6\delta.$$

С учетом того, что (см. утверждение 1)

$$\varphi(y) \ge \varphi(x) + \left\langle \Phi_y(\tilde{x}, y), y - x \right\rangle - \delta \ge \Phi_y(\tilde{x}, y) - 2\delta + \left\langle \Phi_y(\tilde{x}, y), y - x \right\rangle,$$

из определения 1 получаем доказываемое утверждение. ●

Это утверждение, позволяет установить гладкость задачи (2), (3) (но не (5), (6)). Таким образом, для (5), (6) необходимость использования универсального метода для внешней задачи является отражением надежды сходиться быстрее, чем в негладком случае, в то время как для (2), (3) использование универсального метода для внешней задачи является скорее отражением желания настраиваться на правильную константу Липшица градиента. Можно, конечно, пытаться использовать приведенные выше формулы, однако из способа рассуждений (см., например, доказательство утверждения 3) видно, что полученная таким образом константа Липшица может оказаться завышенной.

К сожалению, практическое применение утверждения 3 натыкается на следующие сложности:

1) необходимости отступать от границы множества $Q$ во внутрь на $\sqrt{2\delta/L}$,

2) необходимости рассмотрения ситуации (см. доказательство утверждения 3)

$$\left\| \partial x_i / \partial y_j \right\| = -\Phi_{xx}^{-1}\Phi_{xy},$$

3) необходимости предположения о компактности множества $\bar{Q}$, иначе невозможно будет добиться выполнения условия

$$\max_{z\in Q}\left\langle \Phi_x\left(\tilde{x}, y\right), \tilde{x}-z\right\rangle \le \delta.$$

Сложность 1, как правило, на практике преодолима за счет возможности доопределения функционала задачи с сохранением всех свойств на $\sqrt{2\delta/L}$ окрестность множества $Q$ (заметим, что доопределение часто не требуется, поскольку функционал и так задан "с запасом"). Например, для рассматриваемых нами транспортных приложений с $Q = \left\{ y: \ y \ge \bar{y} \right\}$ это возможно [1 – 4]. Сложность 2 часто вообще не возникает (разве что оговорка о существовании $\Phi_{xx}^{-1}$, впрочем, приведенные выше рассуждения можно провести, сохранив все результаты в идентичном виде, так, что эта оговорка будет не нужна), поскольку $\bar{Q}$ совпадает со всем (двойственным) пространством. А вот сложность 3, действительно, портит дело. К сожалению, простых теоретически обоснованных способов борьбы с этой сложностью мы пока не знаем. Тем не менее, полезно заметить, что в действительности нужно гарантировать выполнение (см. доказательство утверждения 1)

$$\left\langle \Phi_x\left(\tilde{x}\left(y\right), y\right), \tilde{x}\left(y\right)-x\left(y'\right)\right\rangle \le \delta,$$

где точки $y$ и $y'$ близки, поскольку возникают на соседних итерациях внешнего метода. С учетом ожидаемой "близости" $\tilde{x} = \tilde{x}\left(y\right)$ и $x\left(y\right)$, мы можем заменить в этом критерии настоящее множество $\bar{Q}$, которое, как правило, совпадает со всем пространством, на шар конечного радиуса. Более детальные исследования (для задачи (2), (3)) и практические эксперименты показывают, что для выполнения приведенного выше условия достаточно обеспечить для внутреннего итерационного процесса $\left\{ x_k \right\} \to x\left(y\right)$ условия

$$\left\| \Phi_x\left(x_k, y\right) \right\|_2 \left\| x_k \right\|_2 \le \delta/2, \ \left\| \Phi_x\left(x_k, y\right) \right\|_2 \le \delta.$$

Соответствующее $x_k = (\lambda_k, \mu_k)$ порождает нужное $\tilde{x}(y) = x_k$. С учетом специфики рассматриваемой нами задачи (2), (3), имеем следующий критерий (возвращаемся к обозначениям (2), (3))

$$\|Ax(\lambda_k, \mu_k) - b\|_2 \|(\lambda_k, \mu_k)\|_2 \le \delta/2, \|Ax(\lambda_k, \mu_k) - b\|_2 \le \delta,$$

где $x(\lambda_k, \mu_k)$ определяется в формуле (4), а введённая линейная система балансовых уравнений $Ax = b$ – есть общая запись аффинных (транспортных) ограничений:

$$\sum_{j=1}^n x_{ij} = L_i, \sum_{i=1}^n x_{ij} = W_j, \ i, j = 1, ..., n \ .$$

В связи со сказанным выше заметим, что (это следует из оценки (10)) метод балансировки обеспечивает сходимость и по аргументу, что для других методов (без введения регуляризации) решения двойственной задачи, вообще говоря, нельзя гарантировать. Это свойство наряду с линейной скоростью сходимости метода (со скоростью геометрической прогрессии) позволяет надеяться, что выбранный критерий является достаточно точным (точнее не слишком грубым).

Принципиально важно для гладкого случая ($c_{ij}(y)$ – функции с Липшицевым градиентом), как это будет следовать из дальнейших оценок (см. теорему 1), не просто уметь решать двойственную задачу, т.е. находить $(\lambda, \mu)$ так, чтобы была сходимость по аргументу, а делать это так, чтобы сложность решения задачи зависела от точности ее решения логарифмическим образом. Выше мы отмечали, что это имеет место для метода балансировки. Также это имеет место и для быстрых методов, примененных к регуляризованной двойственной задачи. При фиксации параметра регуляризации, исходя из итоговой желаемой точности, быстрые градиентные методы (для сильно выпуклых функций) решают регуляризованную двойственную задачу так, что зависимость сложности от точности ее решения – логарифмическая.

Хочется, чтобы при решении внешней задачи в (2), т.е. задачи

$$\min_{y \in Q} f(y),$$

можно было не задумываться ни о какой гладкости. Если она есть, то метод бы это хорошо учитывал, не требуя знания констант Липшица градиента (это намного более существенно для возможности применять описанный подход к поиску барицентра Вассерштейна вероятностных мер, см. п. 2), если ее нет, то метод также бы работал оптимальным (для негладкого случая) образом. Именно таким свойством и обладает универсальный метод [26], работающий и в концепции неточного оракула [30] (см. определение 1).

Заметим [26], что можно погрузить задачу с гёльдеровым градиентом ($\nu \in [0,1]$)

$$\left\| \nabla f(y') - \nabla f(y) \right\|_* \le L_\nu \left\| y' - y \right\|^\nu$$

(в том числе и негладкую задачу с ограниченной нормой разности субградиентов при $\nu = 0$) в класс гладких задач с оракулом, характеризующимся точностью $\delta$ и

$$L = L_\nu \left[ \frac{L_\nu (1-\nu)}{2\delta(1+\nu)} \right]^{\frac{1-\nu}{1+\nu}}.$$

Это позволяет даже в случае, когда можно рассчитывать только на $\delta$-субградиент[11] (с ограниченной нормой субградиента (разности субградиентов), причем какой именно константой ограниченной, методу знать не обязательно), все равно работать в концепции $(\delta, L)$-оракула.

Итак, у нас есть внешняя задача (2)

$$\min_{y \in Q} f(y),$$

для которой обращение к $(\delta, L)$-оракулу за значением функции и градиента стоит $\sim \ln(\delta^{-1})$. Насколько быстро мы можем решить такую задачу, т.е. при каком $N(\varepsilon)$ можно гарантировать, что

$$f(y_{N(\varepsilon)}) - \min_{y \in Q} f(y) \le \varepsilon ?$$

Ответ можно получить из следующего результата.

**Теорема 1 (см. [20, 26, 30, 34]).** *Существует однопараметрическое семейство универсальных градиентных методов (параметр $p \in [0,1]$), не получающих на вход кроме $p$ больше никаких параметров (в частности, не использующих значения $L_\nu$ и $R$ – "расстояние" от точки старта до решения, априорно не известное), которое приводит к следующей оценке на требуемое число итераций*

$$N_p(\varepsilon) = O\left( \inf_{\nu \in [0,1]} \left( \frac{L_\nu R^{1+\nu}}{\varepsilon} \right)^{\frac{2}{1+2p\nu+\nu}} \right),$$

---

[11] На $\delta$-субградиент всегда можно рассчитывать согласно утверждению 1. Причем, как уже отмечалось раньше, для получения $\delta$-субградиента не нужна сходимость по аргументу для вспомогательной задачи.

*если* $\delta \leq \mathrm{O}\left(\varepsilon / N_p\left(\varepsilon\right)^p\right)$.

Из теоремы 1 можно заключить, что если мы рассчитываем на некоторую гладкость $f\left(y\right)$, то стоит выбирать значение параметра $p = 1$, при этом общее трудозатраты машинного времени будут

$$\mathrm{O}\left(N_1\left(\varepsilon\right)\left(T\ln\left(\varepsilon^{-1}\right)+\tilde{T}\right)\right), \tag{11}$$

где $\tilde{T}$ – время вычисления (суб-)градиента функционала (в основном это вычисления $\left\{\nabla c_{ij}\left(y\right)\right\}_{i,j=1}^{n,n}$ [28, 29]), $T$ – время решения вспомогательной задачи методом балансировки с относительной точностью 1%. Численные эксперименты показывают, что на одном современном ноутбуке при $n \sim 10^2$ время $T \approx 1$ сек. [30], что сопоставимо с временем $\tilde{T}$ для таких $n$ [28].

Выгода, от описанной выше конструкции, по сравнению с обычным способом решения исходной задачи минимизации (2), (3) сразу по совокупности всех переменных (см., например, [4]) заключается в гарантированном не увеличении константы Липшица градиента в оценке необходимого числа итераций (см. утверждение 3) и ожидаемое уменьшение в этой же оценке "расстояния от точки старта до (неизвестного априори) решения. Выгода здесь вполне может достигать одного порядка и более. При этом можно ожидать лишь незначительного увеличения стоимости одной итерации. Причем стоит иметь в виду, что при оптимизации сразу по всем переменным требуется рассчитывать градиента функционала по большему числу переменных, чем в описанном выше подходе, что также играет нам на пользу. В конечном итоге, сокращение числа итераций заметно превалирует над небольшим увеличением стоимости одной итерации.

Что касается задачи (5), (6), то описанный выше подход представляется естественным и не имеющих альтернатив в рассматриваемом классе. Альтернативные методы, с которыми можно сравнивать, мы упоминали по ходу статьи, но все они были предложены из принципиально других подходов. Сравнению (практическому) всех этих методов планируется посвятить отдельную публикацию.

Резюмируем ключевой результат этого раздела (и всей статьи) следующим образом.

Для решения задачи (2) (или (6) или (8)) предлагается использовать универсальный метод из работы [26] (а точнее его модификацию из [30]). Если рассчитываем на гладкость[12] $f(y)$, то полагаем в методе $p = 1$. Если на гладкость рассчитывать не приходится[13], то полагаем $p = 0$. В обоих случаях кроме априорной подсказки относительно параметра $p$, методу больше ничего от нас знать не надо!

## 4. Заключительные замечания

В приложениях часто возникают задачи, имеющие следующий вид (см., например, [35, 36])

$$f(x) = \Phi\big(x, y(x)\big) \to \min_x \qquad (12)$$

при этом $y(x)$ и $\nabla y(x)$ могут быть получены из решения отдельной подзадачи лишь приближенно. Довольно типично, что существует способ, который выдает $\varepsilon$-приближенное решение за время, зависящее от $\varepsilon$, логарифмическим образом $\sim \ln(\varepsilon^{-1})$. В данной работе намечен общий способ решения таких задач. Его наиболее важной отличительной чертой является адаптивность (самонастраиваемость), т.е. методу на вход не надо подавать никаких констант Липшица (более того, метод будет работать и в негладком случае). Метод сам настраивается локально на оптимальную гладкость функции. Это свойство метода делает его привилегированным, поскольку в реальных приложениях, чтобы что-то знать о свойствах $f(x)$ нужно что-то знать о свойствах зависимости $y(x)$, а это часто не доступно по постановке задачи, или, при попытке оценить, приводит к сильно завышенным оценкам.

В задаче (12) важно уметь эффективно пересчитывать значения $y(x)$, а не рассчитывать их каждый раз заново (на каждой итерации внешнего цикла). Поясним сказанное.

---

[12] В этом случае как раз существенна логарифмическая сложность приближенного вычисления градиента и значения функции $f(y)$ от точности, обеспеченная методом балансировки.

[13] В этом случае точность решения вспомогательной задачи расчета $\delta$-субградиента можно завязать на желаемую точность решения задачи (2) $\varepsilon$ (или (6) или (8)) по формуле $\delta = \mathrm{O}(\varepsilon)$ (см. теорему 1 при $\nu = 0$), с константой порядка 1.

Предположим, что мы уже как-то посчитали, скажем, $y(x)$, решив, например, с какой-то точностью соответствующую задачу оптимизации. Тогда для вычисления $y(x + \Delta x)$ (на следующей итерации внешнего цикла) у нас будет хорошее начальное приближение $y(x)$. А, как известно (см., например, [36]), расстояние от точки старта до решения (не в случае сходимости метода со скоростью геометрической прогрессии или быстрей) существенным образом определяет время работы алгоритма. Тем не менее, известные нам приложения (см. [35, 36]) пока как раз всецело соответствуют сходимости процедуры поиска $y(x)$ со скоростью геометрической прогрессии. Связано это с тем, что если расчет $y(x)$ с точностью $\varepsilon$ осуществляется за $O\left(\ln\left(\varepsilon^{-1}\right)\right)$ операций, то для внешней задачи можно выбирать самый быстрый метод (а, стало быть, и самый требовательный к точности), и с точностью до того, что стоит под логарифмом, общая трудоемкость будет рассчитываться по формуле, аналогичной формуле (11). Как правило, такое сочетание оказывается недоминируемым. В частности, для рассматриваемой нами задачи планируется в отдельной статье сравнить описанный здесь подход с возможными альтернативными. Ведь, как известно [20], выпукло-вогнутые седловые задачи, к которым, безусловно, относятся задачи (2), (6), можно по-разному решать, и почему был выбран способ, описанный в п. 3 требует дополнительного сопоставительного анализа, которому и планируется посвятить отдельную статью. Здесь же ограничимся ссылкой на пример 4 и последующий текст из работы [36] и общим тезисом, который пока неплохо подтверждался на практике:

*если есть возможность в задаче оптимизации (в седловой задаче) явно прооптимизировать по части переменных (или как-то эффективно это сделать с хорошей точностью), то, как правило, это и надо сделать, и строить итерационный метод исходя из этого.*

В реальных транспортных приложениях [1, 2, 28, 29] достаточно сложным является расчет $c_{ij}(y)$ и их градиентов (особенно при поиске стохастических равновесий). Тем не менее, эти задачи имеют вполне четкую привязку к решению некоторых задач на графах типа поиска кратчайших путей. Также как и в предыдущем абзаце для внутренней задачи, для внешней задачи, можно не рассчитывать $c_{ij}(y)$ и их градиенты каждый раз заново, а пересчитывать; также можно допускать неточность в их вычислении (и ненулевую вероятность ошибки) надеясь на ускорение (благо метод работает в концепции неточного оракула, природа которой не принципиальна, см. [31]). Также здесь оказываются полезными идеи БАД (быстрого автоматического дифференцирования [37, 38]), которые позволяют практически за тоже время, что занимает вычисление самих функций, вычислять их градиенты.

## Литература